\documentclass[11pt]{article}
 
\usepackage[latin1]{inputenc}
\usepackage{url, amsmath,enumerate,fancyhdr,amssymb, amsthm, 
pstricks, epsf, lscape, ifpdf, graphicx, float}
\usepackage[margin=3.0cm]{geometry}

\renewcommand{\arraystretch}{1.2}

\mathchardef\mhyphen="2D

\newtheorem{Definition}{Definition}
\newtheorem{Example}{Example}
\newtheorem{Proposition}{Proposition}
\newtheorem{Lemma}{Lemma}
\newtheorem{Theorem}{Theorem}
\newtheorem{Corollary}{Corollary}
\newtheorem{Remark}{Remark}
\newtheorem{Assumption}{Assumption}
\newtheorem{Conjecture}{Conjecture}

\newcommand{\lin}{\operatorname{lin}}
\newcommand{\cl}{\operatorname{cl}}
\newcommand{\ri}{\operatorname{ri}}
\newcommand{\dir}{\operatorname{dir}}

\newcommand{\face}{\operatorname{face}}

\newcommand{\Feas}{\operatorname{Feas}}

\newcommand{\ttan}{\operatorname{tan'}}

\newcommand{\beq}{\begin{equation}}
\newcommand{\eeq}{\end{equation}}
\newcommand{\beqa}{\begin{eqnarray}}
\newcommand{\eeqa}{\end{eqnarray}}
\newcommand{\ba}{\begin{array}}
\newcommand{\ena}{\end{array}}
\newcommand{\bac}{\begin{array}{ccccccccccc}}
\newcommand{\eac}{\end{array}}
\newcommand{\bprop}{\begin{Proposition}}
\newcommand{\eprop}{\end{Proposition}}
\newcommand{\beqast}{\begin{eqnarray*}}
\newcommand{\eeqast}{\end{eqnarray*}}
\newcommand{\benum}{\begin{enumerate}}
\newcommand{\eenum}{\end{enumerate}}
\newcommand{\bit}{\begin{itemize}}
\newcommand{\eit}{\end{itemize}}
\newcommand{\bth}{\begin{Theorem}}
\newcommand{\enth}{\end{Theorem}}
\newcommand{\ble}{\begin{Lemma}}
\newcommand{\ele}{\end{Lemma}}
\newcommand{\bex}{\begin{Example}}
\newcommand{\eex}{\end{Example}}
\newcommand{\bcor}{\begin{Corollary}}
\newcommand{\ecor}{\end{Corollary}}
\newcommand{\brem}{\begin{Remark}}
\newcommand{\erem}{\end{Remark}}
\newcommand{\bass}{\begin{Assumption}}
\newcommand{\eass}{\end{Assumption}}
\newcommand{\bsmx}{\begin{small} \begin{pmatrix}}
\newcommand{\esmx}{\end{pmatrix} \end{small}}
\newcommand{\bpx}{\begin{pmatrix}}
\newcommand{\epx}{\end{pmatrix}}
\newcommand{\bbx}{\begin{bmatrix}}
\newcommand{\ebx}{\end{bmatrix}}
\newcommand{\bdef}{\begin{Definition}} 
\newcommand{\bconj}{\begin{Conjecture}}
\newcommand{\econj}{\end{Conjecture}}
\newcommand{\bscr}{\begin{scriptsize}}
\newcommand{\escr}{\end{scriptsize}}
\newcommand{\bfoot}{\begin{footnotesize}}
\newcommand{\efoot}{\end{footnotesize}}
\newcommand{\bsmall}{\begin{small}}
\newcommand{\esmall}{\end{small}}

\newcommand{\commentout}[1]{}
\newcommand{\co}[1]{}

\newcommand{\nin}{\noindent}
\newcommand{\ti}{\times}
 
\newcommand{\pf}[1]{\vspace{.35cm} \nin {\bf Proof {#1} }}

\newcommand{\sym}[1]{{\cal S}^{#1}}

\newcommand{\psd}[1]{{\cal S}_+^{#1}}
\newcommand{\rad}[1]{\mathbb{R}^{#1}}
\newcommand{\radp}[1]{\mathbb{R}^{#1}_+}

\newcommand{\eref}[1]{(\ref{#1})}

\newcommand{\la}{\langle}
\newcommand{\ra}{\rangle}

\newcommand{\tri}{\triangle}

\newcommand\conelpspace{\hspace{.9cm}}
\newcommand{\redf}{\mathit{R \mhyphen P}}
\newcommand{\redfd}{\mathit{R \mhyphen D}}
\newcommand{\fra}{\operatorname{FRS}}

\begin{document}
\title{\bf  Strong duality in conic linear programming: facial reduction and extended duals} 
\author{G\'{a}bor Pataki  \\ {\bf gabor@unc.edu} \vspace{.5cm} \\ Department of Statistics and Operations Research \\
University of North Carolina at Chapel Hill }

\maketitle


\renewcommand{\qedsymbol}{$\blacksquare$}

\pagestyle{myheadings}
\thispagestyle{plain}
\markboth{}{Strong duality in conic linear programming: facial reduction and extended duals}

\begin{abstract}
The facial reduction algorithm of Borwein and Wolkowicz and 
the extended dual of Ramana provide a strong dual for the conic linear program 
\begin{equation} \tag{$P$} \label{P}
\sup \, \{ \, \langle c, x \rangle \, | \, Ax  \leq_K  b \, \}
\end{equation}
in the absence of any constraint qualification. 
The facial reduction algorithm solves a sequence of auxiliary optimization problems to obtain such a dual.
Ramana's dual is applicable when (\ref{P}) is 
a semidefinite program (SDP) and is an explicit SDP itself.
Ramana, Tun\c{c}el, and Wolkowicz showed that these approaches are closely related;
in particular, they proved the correctness of Ramana's dual using certificates from 
a facial reduction algorithm.
Here we give a simple and self-contained 
exposition of facial reduction, of extended duals, and generalize Ramana's dual:
\begin{itemize}
\item we state a simple facial reduction algorithm and prove its correctness; and 
\item building on this algorithm we construct a family of extended duals 
when $K$ is a {\em nice} cone. 
This class of cones includes the semidefinite cone and other important cones.
\end{itemize}
\end{abstract}

{\em Dedicated to Jonathan Borwein on the occasion of his 60th birthday}



\tableofcontents

\section{Introduction}

Conic linear programs generalize ordinary linear programming, as they require
membership in a closed convex cone in place of the usual nonnegativity constraint.
Conic LPs 
share some of the duality theory of linear optimization: weak duality always holds 
in a primal-dual pair, and assuming a suitable constraint qualification (CQ), their 
objective values agree, and are attained. 

When a CQ is lacking and the 
underlying cone is not polyhedral, pathological phenomena can occur: 
nonattainment of the optimal values, positive gaps, and infeasibility of the dual even when the 
primal is bounded. All these pathologies appear in semidefinite programs,
second order cone programs, and $p$-order conic programs, arguably the most important and useful 
classes of conic LPs 
(\cite{BonnShap:00, BentalNem:01, Todd:00, AliGold:03, 
SaigVandWolk:00, BoydVand:04, Guler:10, AndRoosTerlaky:02}).

\paragraph[fac-par]{\bf Facial reduction and extended duals} 
Here we study two fundamental approaches to duality in conic 
linear programs that work without assuming any CQ. The first approach is the 
facial reduction algorithm (FRA) of Borwein and Wolkowicz \cite{BorWolk:81, BorWolk:81B}, 
which constructs a so-called 
minimal cone of a conic linear system. 
Using this minimal cone one can always ensure strong duality in a primal-dual pair of conic LPs.

The second approach is Ramana's extended dual
for semidefinite programs \cite{Ramana:97}.
(Ramana named his dual an Extended Lagrange-Slater Dual, or ELSD dual. 
We use the shorter name for simplicity.) 
The extended dual is an explicit semidefinite program with 
a fairly large  number (but polynomially many) variables and constraints.
It has the following desirable properties: it is feasible if and only if the 
primal problem is bounded; and when these equivalent statements hold, it 
has the same value as 
the primal and attains it. 

Though these approaches at first sight look quite different, Ramana, Tun\c{c}el,
and Wolkowicz in \cite{RaTuWo:97} showed that 
they are closely related: 
in the case of semidefinite programming, they proved the  correctness of Ramana's dual 
using certificates from the algorithm of \cite{BorWolk:81, BorWolk:81B}.

The goal of our paper is to give a simple and self-contained
exposition of facial reduction, of extended duals, study their connection, 
and give simple proofs of generalizations of Ramana's dual. 
We use ideas from the paper of Ramana, Tun\c{c}el, and Wolkowicz \cite{RaTuWo:97}, although
our development is different. 
We state a facial reduction algorithm 
and prove its correctness using only elementary results from the duality theory of conic LPs, and 
convex analysis. We build on this algorithm and 
generalize Ramana's dual: we construct a family of extended duals
for (\ref{P}) when $K$ is a {\em nice} cone. This class of cones includes the semidefinite 
cone, and other important cones, as $p$-order, in particular, second order cones. 

Next we present our framework in more detail. 
A conic linear program can be stated as 
\begin{equation} \tag{$P$} \label{p}
\begin{array}{lrccrcr}                                                                                                   & \sup  & \langle c, x \rangle   \\
      & s.t.  & Ax \leq_K b,   \\ 
\end{array}                                                                                           
\end{equation}
where $A: X \rightarrow Y$ is a linear map between 
finite dimensional Euclidean spaces 
$X$ and $Y, \,$ and $c \in X, \, b \in Y. \,$ 
The set $K \subseteq Y$ is a closed, convex cone, and we write 
\mbox{$Ax \leq_K b \,$} to mean $b - Ax \in K.$ 
We naturally associate a dual program with (\ref{p}). 
Letting $A^*$ be the adjoint operator of $A, \,$ and $K^*$ the dual cone of $K, \,$ i.e., 
$$
K^* \, = \, \{ \, y \, | \, \langle y, x \rangle \geq 0 \, \forall x \in K \, \},
$$
the dual problem is 
\begin{equation}  \tag{$D$} \label{d}
\begin{array}{rrl} 
       & \inf        &   \langle b, y \rangle     \\                   
       & s.t.        &  y \geq_{K^*} 0 \\
       &             & A^*y  = c.   
\end{array}
\end{equation}
When (\ref{p}) is feasible,  
we say that {\em strong duality holds between (\ref{p}) and (\ref{d})} if the following conditions 
are satisfied:
\begin{itemize}
\item problem (\ref{p}) is bounded, if and only if (\ref{d}) is feasible; and 
\item when these equivalent conditions hold, the optimal values of (\ref{p}) and (\ref{d}) agree and the 
latter is attained.
\end{itemize}
We say that (\ref{p}) is {\em strictly feasible}, or {\em satisfies Slater's condition}, 
if there is an $x \in X$ such that $b - Ax$ is in the relative interior of $K. \,$ 
When (\ref{p}) is strictly feasible, it is well-known that strong duality holds between 
(\ref{p}) and (\ref{d}). 

The facial reduction algorithm (FRA) of Borwein and Wolkowicz constructs a suitable face of 
$K, \,$ called the {\em minimal cone of (\ref{P})}, which we here denote by $F_{\min}.$ 
The minimal cone has two important properties:
\begin{itemize}
\item The feasible set of (\ref{P}) remains the same if we replace its constraint set
by 
$$
Ax \leq_{F_{\min}} b.
$$
\item The new constraint set satisfies Slater's condition. 
\end{itemize}
Thus, if we also replace 
$K^*$ by $F_{\min}^*$ in (\ref{d}), strong duality holds in the 
new primal-dual pair. 
The algorithm in \cite{BorWolk:81, BorWolk:81B} 
constructs  a decreasing chain of faces starting with $K$ and ending with 
$F_{\min},$ in each step solving a pair of auxiliary conic linear programs. 

\paragraph[contr-par]{\bf Contributions of the paper} 
We first state a simplified FRA and prove its correctness. 
Building on this algorithm, and assuming that 
cone $K$ is {\em nice, \,} i.e., the set $K^* + F^\perp$ is closed for all $F$ faces 
of $K \,$ we show that the dual of the minimal cone has a
representation
\begin{equation} \label{fmin*repr}
\begin{array}{rcl}
F_{\min}^* & = & \{ \, u_{\ell+1} + v_{\ell+1} \, : \, \\
           &   & \, (u_0, v_0) = (0, 0) \\
           &   &    \,    (A, b)^*(u_i + v_i) = 0, \, i=1, \dots, \ell, \\
           &   &    \,   (u_{i}, v_{i}) \in K^* \times \tan(u_0 + \dots + u_{i-1}, K^*), \,i=1, \dots, \ell+1 \, \},
\end{array}
\end{equation}
where $\tan(u, K^*)$ denotes the tangent space of the cone $K^*$ at $u \in K^*$ 
and $\ell$ is a suitable integer.
Plugging this expression 
for $F_{\min}^*$ in place of $K^*$ in (\ref{d}) we obtain a dual with the properties of Ramana's dual.
We show the correctness of several representations of $F_{\min}^*,$ 
each leading to a different extended dual. We note that the results of \cite{RaTuWo:97} 
already imply that such a representation is possible, but this is not stated there explicitly.

The cone of positive semidefinite matrices is nice
(and also self-dual), so in this case the representation of (\ref{fmin*repr}) is valid.
In this case we can 
also translate the tangent space 
constraint into an explicit semidefinite constraint and recover variants of 
Ramana's dual.

We attempted to simplify our treatment of the subject as much as possible:
as background we use only the fact that strong duality holds in a primal-dual pair of conic LPs, 
when the primal is strictly feasible and some elementary facts in convex analysis. 

\paragraph[lit-rev]{\bf Literature review} Borwein and Wolkowicz originally presented 
their facial reduction algorithm 
in \cite{BorWolk:81} and \cite{BorWolk:81B}. Their algorithm works for a potentially nonlinear 
conic system of the form $\{ \, x \, | \, g(x) \in K \, \}.$
Luo, Sturm, and Zhang in 
\cite{LuoSturmZhang:97} studied a so-called {\em conic expansion method} which finds 
a sequence of increasing sets starting with $K^*$ and ending with $F_{\min}^*:$ 
thus their algorithm can be viewed as a dual variant of facial reduction. 
Their paper also contains an exposition of facial reduction and Ramana's dual.
Sturm in \cite{Sturm:00} introduced an 
interesting and  novel application of facial reduction: 
deriving error bounds for semidefinite systems that lack a strictly feasible solution. 
Luo and Sturm in \cite{LuoSturm:00} generalized this approach to mixed semidefinite
and second order conic systems.
Lewis in \cite{Lewis:94} used facial reduction to 
derive duality results without a CQ assumption in partially finite convex programming.
Tun\c{c}el in his recent book \cite{Tuncel:11} constructed an SDP instance with 
$n$ by $n$ semidefinite matrices that requires $n-1$ iterations of the facial reduction algorithm to find 
the minimal cone, and thus showed that the theoretical worst case is essentially attainable.

Waki and Muramatsu in \cite{WakiMura:12} also described an FRA, rigorously 
showed its equivalence to the conic expansion 
approach of Luo et al, and presented computational results on semidefinite programs. 
A preliminary version of this paper appeared in \cite{Pataki:00B}. P\'olik and Terlaky 
in \cite{PolikTerlaky:09} used the results of \cite{Pataki:00B} to construct
strong duals for conic LPs over homogeneous cones. Wang et al in \cite{WangXiuLuo:11} 
presented a facial reduction algorithm for nonsymmetric semidefinite least squares problems. 

Tun\c{c}el and Wolkowicz in 
\cite{TunWolk:12} described a connection between the lack of strict complementarity in the 
homogeneous primal and dual systems, and positive duality gaps in SDPs: in particular, they proved 
that when strict complementarity in the homogeneous problems fails in a certain minimal sense, 
one can generate instances with an arbitrary positive duality gap.
Waki in \cite{Waki:12} showed how to systematically 
find SDP instances that are {\em weakly infeasible}, i.e., 
infeasible without a Farkas' lemma certificate. 
Cheung et al in \cite{CheWolkSchurr:12} developed a relaxed version of a facial reduction algorithm,
in which one can allow an error in the solution of the auxiliary conic LPs, 
and applied their method to SDPs, in particular, to instances generated according to 
the results of \cite{TunWolk:12}.

Nice cones appear in other areas of optimization as well. 
In \cite{Pataki:07} we studied the question of when the linear image of a closed convex cone 
is closed and described necessary and sufficient conditions. These lead to 
a  particularly simple and exact characterization when the dual of the cone in question is nice.
We call a conic linear system {\em well behaved} if for all objective functions the 
resulting conic linear program has strong duality with its dual and {\em badly behaved}, otherwise.
In related work, \cite{Pataki:10}, 
we described characterizations of well- and badly behaved conic linear systems. These become 
particularly  simple when the underlying cone is nice, and yield 
combinatorial type characterizations for semidefinite and second order conic systems.

Chua and Tun\c{c}el in \cite{ChuaTuncel:08} showed that if a cone $K$ is nice, then 
so is its intersection with a linear subspace. 
Thus, all homogeneous cones are nice, since
they arise as the slice of a semidefinite cone with a suitable subspace, as proven 
independently by Chua in \cite{Chua:03} and by Faybusovich in \cite{Faybu:02}. 
In \cite{ChuaTuncel:08} the authors also proved that 
the preimage of a nice cone under a linear map is also nice and 
in \cite{Pataki:12} we pointed out that this result implies 
that the intersection of nice cones is also nice. In \cite{Pataki:12} 
we gave several characterizations of nice cones and proved that 
they must be facially exposed; facial exposedness with a 
mild additional condition implies niceness;  
and conjectured that facially exposed and nice cones are actually the same class of cones. 
However, Roshchina disproved this conjecture \cite{Vera:13}.

Most articles on strong duality deal with instances with a fixed right hand side.
Schurr et al in \cite{Schurretal:07} obtained characterizations of 
{\em universal duality}, i.e., of the situation when strong duality holds
for all right hand sides, and objective functions.

Klep and Schweighofer in \cite{KlepSchw:12} derived a strong dual for semidefinite programs
that also works without assuming any constraint qualification.
Their dual resembles Ramana's dual. However, it is based on 
concepts from algebraic geometry, whereas all other references known to us use 
convex analysis.

Recently Gouveia et al in \cite{GouveiaParriloThomas:12} studied the following fundamental question:
can a convex set be represented as the projection of an affine slice of a suitable closed, convex cone?
They gave necessary and sufficent conditions for such a {\em lift} to exist and showed that 
some known lifts from the literature are in the lowest dimension possible. The representation 
of (\ref{fmin*repr}) is related in spirit, as we also represent the set
$F_{\min}^*$ as the projection of a conic linear system in a higher dimensional 
space. 

\paragraph[org-par]{\bf Organization of the paper and guide to the reader} 
In Section \ref{sect-prelim} we fix notation,  review 
preliminaries, and present two motivating examples. 
The reader familiar with convex analysis can 
skip the first part of this section and go directly to the examples.
In Section \ref{fra-section}  we present a simple facial reduction algorithm, prove its
correctness, and show how $F_{\rm min}^*$ can be written as the projection of a nonlinear
conic system in a higher dimensional space. 

Assuming that $K$ is nice, in Section \ref{section-nice} we 
arrive at the representation in (\ref{fmin*repr}), i.e., show that 
$F_{\rm min}^*$ is the projection of a conic {\em linear} system,
and derive an extended dual for conic LPs over nice cones. 
Here we obtain our first Ramana-type dual for semidefinite programs
which is an explicit SDP itself, but somewhat different from the dual proposed in \cite{Ramana:97}.

In Section \ref{section-variants} we describe variants of the representation in (\ref{fmin*repr}), 
of extended duals, and show how we can exactly obtain Ramana's dual. 
In Section \ref{repr-sect} we show that 
the minimal cone  $F_{\min} \,$ itself also has a representation similar to 
the representation of $F_{\min}^*$ in (\ref{fmin*repr}) and discuss some open questions.

The paper is organized to arrive at an explicit Ramana-type dual for SDP as quickly as possible.
Thus, if a reader is interested in  only the derivation of such a dual, it suffices for him/her
to read only Sections \ref{sect-prelim}, \ref{fra-section} and \ref{section-nice}.

\section{Preliminaries}
\label{sect-prelim} 
\paragraph[mat-par]{Matrices and vectors} 
We denote operators by capital letters and matrices (when they are considered as elements 
of a Euclidean space and not as operators) and vectors by lower case letters. 
The $i$th component of vector $x$ is denoted by $x_i \,$ and the $(i,j)$th component of matrix 
$z$ by $z_{ij}.$ 
We distinguish vectors and matrices of similar type with lower indices, i.e., writing
$x_1, x_2, \dots$  The $j$th component of vector $x_i$ is denoted by
$x_{i,j}. \,$ This notation 
is somewhat ambiguous, as $x_i$ may denote a vector, or the $i$th component of
the vector $x, \,$ but the context will make it clear which one is meant. 

\paragraph[convex-par]{Convex sets} 
For a set $C$ we write $\cl C$ for its closure, $\lin C$ for its linear span,
and $C^\perp$ for the orthogonal complement of its linear span. 
For a convex set $C$ we denote its relative interior by $\ri C.$
For a one-element set $\{ x \}$ we abbreviate $\{ x \}^\perp$ by $x^\perp.$  
The open line-segment between points $x_1$ and $x_2$ is denoted by $(x_1,x_2)$. 

For a convex set $C, \,$ and an $F, \,$ a convex subset of $C \,$ 
we say that $F$ is a face of $C \,$ if
$x_1, x_2 \in C$ and $(x_1, x_2) \cap F \neq \emptyset$ implies that 
$x_1$ and $x_2$ are both in $F.$
For $x \in C$ there is a unique minimal face of $C$ that contains $x, \,$ i.e., the face that contains 
$x$ in its relative interior: we denote this face by  $\face(  x , C).$ 
For $x \in C$ we define the set of feasible directions and the tangent space at $x$ in $C$ as 
\begin{align*}
\dir(x,C) & = \, \{\, y \, | \, x + t y \in C \; {\rm for \; some \;} t > 0 \, \}, \\ 
\tan(x,C) & = \, \cl \dir(x,C) \cap  - \cl \dir(x,C).  
\end{align*}

\paragraph[cone-par]{Cones}  We say that a set $K$ is a  cone, if  $\lambda x \in K$ holds for all 
$x \in K \,$ and $\lambda \geq 0,$ and define the dual of cone $K$ as
\begin{align*}
K^*  & = \, \{\, z \, | \, \la z, x \ra \geq 0 \; {\rm for \; all \;} x \in K \, \}. 
\end{align*}
For an $F$ face of  a closed convex cone $K$ and $x \in \ri F$ 
the  complementary, or conjugate face of $F$ is defined alternatively as 
(the equivalence is straightforward)
\begin{align*}
F^\tri    & = \, K^* \cap F^\perp \, = \,  K^* \cap x^\perp. 
\end{align*}
We define the complementary face of a face $G$ of $K^*$ 
analogously and denote it by $G^\triangle$. 
A closed convex cone $K$ is facially exposed, 
i.e., all faces of $K$ arise as the intersection of $K$ with a 
supporting hyperplane iff for all $F$ faces of $K$ we have $(F^{\tri})^\tri = F.$
For brevity  we write $F^{\tri *}$  for $(F^{\tri})^*$ and $F^{\tri \perp}$ for $(F^{\tri})^\perp$. 

For a closed convex cone $K $ and $x \in K$ we have 
\beqa 
\label{tanxK}
\tan(x, K)         & = & \face(x, K)^{\tri \perp},
\eeqa
as shown in \cite[Lemma 1]{Pataki:10}. 

\paragraph[psd-para]{The semidefinite cone} 
We denote the space of $n$ by $n$ symmetric and the cone of 
$n$ by $n$ symmetric, positive semidefinite  matrices by $\sym{n},$ and $\psd{n}$, respectively. 
The space $\sym{n}$ is equipped with the inner  product 
\beqast
\la x , z \ra & := & \sum_{i,j=1}^n x_{ij} z_{ij}, 
\eeqast
and $\psd{n}$ is self-dual with respect to it. 
For  $y \in \sym{n}$ we write $y \succeq 0$ to denote that $y$ is positive semidefinite.
Using a rotation $v^T (.) v$ by a full-rank matrix $v \,$ 
any face of $\psd{n}$ and its conjugate face can be brought to the form 
\beq \label{psd-face}
F \, = \, \Biggl\{ \bpx x & 0 \\ 0 & 0 \epx \, | \,  x \in \psd{r} \, \Biggr\}, \, F^\tri \, = \, \Biggl\{ \bpx 0 & 0 \\ 0 & y \epx \, | \,  y \in \psd{n-r} \, \Biggr\},
\eeq
where $r$ is a nonnegative integer. 

For a face of this form and related sets  we use the shorthand
\beq \label{f-oplus}
F = \begin{pmatrix}  \oplus \! &  0  \\
                    0   \! &  0 
     \end{pmatrix}, \, F^\tri = \begin{pmatrix}  0 &  0  \\
                    0   \! &  \oplus \!  
     \end{pmatrix}, \, F^{\tri \perp} = \begin{pmatrix}  \ti &  \ti  \\
                    \ti   \! &  0   
     \end{pmatrix}, \, 
\eeq
when the size of the partition is clear from the context. 
The $\oplus$ sign denotes a positive semidefinite submatrix and the sign
 $\ti$ stands for a submatrix with arbitrary elements. 

For an $x$ positive semidefinite matrix 
we collect some expressions for $\tan(x, \psd{n})$ below: these play 
an important role when constructing explicit duals for semidefinite programs. 
The second part of Proposition \ref{tan-prop} is based on Lemma 1 in \cite{RaTuWo:97}.
\bprop \label{tan-prop}
The following statements hold.
\benum
\item \label{tan-prop-1} Suppose $x \in \psd{n}$ is \co{a positive semidefinite matrix} of the form
\beq \label{xir}
x \, = \, \bpx I_r & 0 \\ 0 & 0 \epx,
\eeq
and $F = \face(x, \psd{n}).$ 
Then $F, \, F^\tri,$ and $F^{\tri \perp}$ are as displayed in equation 
\eref{f-oplus}, with the upper left block 
$r$ by $r$, and 
\beq \label{tanir}
\tan(x, \psd{n}) \, = \, F^{\tri \perp} \, = \, \bpx  \ti &  \ti  \\
                    \ti   \! &  0   
     \epx. 
\eeq
\item \label{tan-prop-2} For an arbitrary $x \in \psd{n}$ we have 
\beq \label{tan-x}
\tan(x, \psd{n})  = \left\{ \, w + w^T \Bigm |  \, \bpx x    & w \\
                                                        w^T  & \beta I
		         \end{pmatrix} \,  \succeq 0   \,\, \text{for some} \;   
\beta  \in  \rad{} \right\} .
\eeq
\eenum
\eprop
\pf{of \eref{tan-prop-1}} This statement is straightforward from the form of $x$ and the expression 
for the tangent space given in \eref{tanxK} with $K = \psd{n}.$ 

\pf{of \eref{tan-prop-2}} 
If  $x$ is of the form as in equation \eref{xir}, then our claim 
follows from part \eref{tan-prop-1}. 

Suppose now that $x \in \psd{n}$ is arbitrary 
and let $q$ be a matrix 
of suitably scaled eigenvectors of $x$ with eigenvectors corresponding to nonzero eigenvalues coming first.
Let us write $T(x)$ for the set on the right hand side of equation \eref{tan-x}. 
Then one easily checks $\tan(q^T x q, \psd{n}) = q^T \tan(x, \psd{n})q$ and 
$T(q^T x q) = q^T T(x) q, \,$ so this case reduces to the previous case.
\qed

\paragraph[sdp-par]{Conic LPs}
An ordinary linear program is clearly a special case of \eref{p}.
If we choose $X = \rad{m}, \, Y = \sym{n}, \,$ and $K = \psd{n}, \,$ then problem \eref{p} 
becomes a semidefinite program (SDP). 
Since $K$ is self-dual, the dual problem \eref{d} is also an SDP.
The operator $A$ and its adjoint are defined via 
symmetric matrices $a_1, \dots, a_m$ as 
$$
Ax \, = \, \sum_{i=1}^m x_i a_i \,\, {\rm and} \,\, A^*y = (\la a_1, y \ra, \dots, \la a_m, y \ra)^T.
$$
We use the operator $\Feas()$ to denote the feasible set of a conic system.

\paragraph[minc-par]{The minimal cone}
Let us choose $x \in \ri \Feas(P).$ 
We define the minimal cone of \eref{P}  as the unique face of $K$ that contains 
$b - Ax$ in its relative interior and denote this face by $F_{\min}.$ 

For an arbitrary $y \in \Feas(P)$ there is $z \in \Feas(P)$ such that $x \in (y,z)$.
Hence $b - Ax \in (b - Ay,b - Az), \,$ so $b - Ay$ and $b - Az$ are in $F_{\min}, \,$ and 
\eref{P} is equivalent to 
\begin{align*}
Ax & \leq_{F_{\min}}  b,
\end{align*}
and this constraint system satisfies Slater's condition.

\paragraph[nice-par]{Nice cones} 
We say that a closed convex cone $K$ is  nice if 
\beq
\nonumber 
K^* + E^\perp   \; {\rm is \; closed } \; {\rm for \; all \; } E \; {\rm faces \; of \;} K.
\eeq
Most cones appearing in the optimization literature, such as 
polyhedral, semidefinite, $p$-order, in particular 
second order cones are nice: see e.g. \cite{BorWolk:81, BorWolk:81B, Pataki:07}. Furthermore, geometric and dual 
geometric cones are nice as well \cite{Glineur:01}. 

\bex \label{lin-ex} {\rm In the linear inequality system
\beq \label{linsystem}
\begin{pmatrix}  1 & 0  & 0 \\
	         0 & -1 & 1 \\
	         0 & 1  & 0 \\
	         0 & 0  & -1 \\ 
	         0 & 0  &  1 \\
\end{pmatrix} 
\begin{pmatrix}  x_1 \\
	         x_2 \\
	         x_3
\end{pmatrix} 
\leq 
\begin{pmatrix}  0 \\
	          0 \\
	          0 \\
	          0 \\
	          0 
\end{pmatrix} 
\end{equation}
all feasible solutions satisfy the last four inequalities at equality, and for, say,
 $x = (-1, 0, 0)^T$ the first inequality is strict. So the minimal cone of this system is
$$
F_{\rm min} = \radp{1} \ti \{ 0 \}^4.
$$ 
In linear programs strong duality holds even without strict feasibility, so  this example
illustrates only the concept of 
the minimal cone. 

} \qed
\eex

\bex \label{ex1} {\rm 
In the semidefinite program
\beq \label{ex1-p}
\ba{rl}
\sup &  x_1 \\
s.t. & x_1 \bpx 0 & 1 & 0 \\ 1 & 0 & 0 \\ 0 & 0 & 0 \epx + x_2 \bpx 0 & 0 & 1 \\ 0 & 1 & 0 \\ 1 & 0 & 0 \epx \preceq \bpx 1 & 0 & 0 \\ 0 & 0 & 0 \\ 0 & 0 & 0 \epx
\ena
\eeq
a feasible positive semidefinite 
slack $z$ must have all entries equal to zero, except for $z_{11},$ and 
there is a feasible slack with $z_{11}>0.$
So the minimal cone and its dual are
\beq \label{fminex1}
F_{\min} \, = \, \bpx \oplus & 0 & 0 \\
                        0    & 0 & 0 \\
                        0    & 0 & 0 \epx, \, 
F_{\min}^* \, = \, \bpx \oplus & \ti & \ti \\
                        \ti    & \ti & \ti \\
                        \ti    & \ti & \ti \epx. \, 
\eeq
The optimal value of \eref{ex1-p} is clearly zero. Writing $y$ for the dual matrix, 
the dual program is equivalent to 
\beq \label{ex1-d} 
\ba{rl}
\inf & y_{11} \\
s.t. & \bpx y_{11} & 1/2        & - y_{22}/2 \\
             1/2     & y_{22}   & y_{23} \\
        - y_{22}/2     &  y_{23}  & y_{33} \epx \succeq 0.
\ena
\eeq
The dual has an unattained $0$ minimum: $y_{11}$ can be an arbitarily small positive number, 
at the cost of making $y_{22}$ and in turn $y_{33}$ large, however, 
$y_{11}$ cannot be $0,$ as $y_{12}$ is $1/2.$ 

Suppose that in \eref{ex1-d}  we replace the constraint $y \succeq 0$ by 
$y \in F_{\min}^*. \,$ Then we can set $y_{11}$ to zero, 
so with this modification the dual attains. } \qed
\eex 

We will return to these examples later to illustrate our facial reduction algorithm and 
extended duals.

We assume throughout the paper that \eref{P} is feasible. It is possible to remove this assumption
and modify the facial reduction algorithm of section \ref{fra-section} to either prove the infeasibility of 
\eref{P}, or to find the minimal cone in finitely many steps; such an FRA was described by 
Waki and Muramatsu in \cite{WakiMura:12}.

\section{A simple facial reduction algorithm}

\label{fra-section} 
We now state a simple facial reduction algorithm that is applicable 
when $K$ is an arbitrary closed convex cone.
We prove its correctness and illustrate it on Examples \ref{lin-ex} and \ref{ex1}. 

Let us recall that $F_{\min}$ denotes the minimal cone of (\ref{p}) and for brevity define the subspace $L$ as
\begin{equation} \label{def-L}
L = {\cal N}( (A, b)^*).
\end{equation}
\ble 
\label{red-lemma}
Suppose that an $F$ face of $K$ satisfies $F_{\min} \subseteq F.$ 
Then the following hold:
\begin{enumerate} 
\item \label{red-lemma-1} For all $y \in F^* \cap L$ we have 
\begin{equation} \label{FminF}
F_{\min} \subseteq F \cap y^\perp \subseteq F.
\end{equation}
\item \label{red-lemma-2} 
There exists $y \in F^* \cap L$ such that the second containment in 
(\ref{FminF}) is strict, iff $F_{\min} \neq F. \,$ 
We can find such a $y, \,$ or prove $F = F_{\min}$ by solving a pair of auxiliary conic linear programs.
\end{enumerate}
\ele
\pf{of \eref{red-lemma-1}} To prove (\ref{red-lemma-1}) suppose that $x$ is feasible for $(P) \,$ and let 
$y \in F^* \cap L.$ Then 
$b - Ax \in F_{\min} \subseteq F, \,$ hence $\langle b - Ax, y \rangle = 0, \,$ which implies the first 
containment; the second is obvious. 

\pf{of \eref{red-lemma-2}} The ``only if'' part of the statement 
is obvious. To see the ``if'' part, let us fix $f \in \ri F, \,$ 
and consider the primal-dual pair of conic linear 
programs that we call reducing conic LPs below:
\begin{center}                                                                                                 $$     \begin{array}{lrccrcr}                                                                                      
       &   \sup  & t                 &  \conelpspace & \inf    &  \langle b, y \rangle    &   \\                   
(\redf)    &  s.t.   & Ax + ft \leq_F b  &  \conelpspace &  s.t.   &  y \geq_{F^*} 0  & (\redfd) \\
       &         &                   &  \conelpspace &         &   A^* y = 0      &     \\
       &         &                   & \conelpspace &          &   \langle f, y \rangle = 1.      &    
\end{array}                                                                                              
$$   
\end{center}                                                                                                      
First let us note 
$$
\begin{array}{rcl} 
F_{\min} = F & \Leftrightarrow & \exists x \; \text{s.t.} \, b - Ax \in \ri F \\
             & \Leftrightarrow & \exists x \; \text{and} \, t > 0 \, \text{s.t.} \, b - Ax - ft \in F.
\end{array}
$$
Here in the first equivalence the direction $\Rightarrow$ is obvious from the definition of the minimal cone.
To see the direction $\Leftarrow$ assume $b - Ax \in \ri F.$ Then $\ri F \cap F_{\min} \neq \emptyset \,$ 
and  $F_{\min}$ is a face of $K, \,$ so Theorem 18.1 in 
\cite{Rockafellar:70} implies
$F \subseteq F_{\min}, \,$ and the reverse  containment is already given.
The second equivalence is obvious.

Therefore, $F_{\min} \neq F$ iff the optimal value of $(\redf)$ is $0.$ 
Note that $(\redf)$ 
is strictly feasible, with some $x$ such that $b - Ax \in F, \,$ and some $t < 0. \,$

Hence $F_{\min} \neq F$ iff $(\redfd)$ has optimal value $0 \,$ and attains it, i.e., iff there is
$y \in F^* \cap L$ with $ \langle f, y \rangle =1.$ 
Such a $y$ clearly must satisfy $F \cap y^\perp \subsetneq F.$ 
\qed

Based on Lemma \ref{red-lemma} we now state a simple facial reduction algorithm in Figure 
\ref{fig:fr-alg}.
\vspace{.35cm}
\begin{figure}[ht]
\framebox[5.25in]{\parbox{5.05in}{ 
{\sc Facial Reduction Algorithm}
\begin{tabbing}{}
**************\=******\=*******\=***\= \hspace{2.5in} \=  \kill  
 {\bf Initialization:} \> Let $y_0 = 0, \, F_0 = K, \, i=1.$ \\
 {\bf repeat }  \\
 \> Choose \= $y_{i} \in L \cap F_{i-1}^*.$                 \\
 \>        \= Let $F_i = F_{i-1} \cap y_i^\perp.$  \\
 \> Let  \= $i = i+1.$ \\
 {\bf end repeat}    \\
\end{tabbing}
}}
\caption{The facial reduction algorithm} 
\label{fig:fr-alg}
\end{figure}
\vspace{.35cm}

The algorithm of Figure \ref{fig:fr-alg} may not terminate in general, 
as it allows the choice of 
a $y_i$ in iteration $i$ such that $F_{i} = F_{i-1};$ it even allows $y_i = 0$ for all $i$.
Based on this general algorithm, however, it will be convenient to 
construct a representation of $F_{\min}^*.$ 

We call an iteration of the FRA reducing, if the $y_{i}$ vector found therein satisfies 
$F_{i} \subsetneq F_{i-1};$ 
we can make sure that an iteration is reducing, or that we have found the minimal cone 
by solving the pair of conic linear programs $(\redf)\mhyphen(\redfd).$ 
It is clear that after a sufficient number of reducing iterations
the algorithm terminates.

Let us define the quantities 
\begin{equation} \label{def-ell}
\begin{array}{rcl}
\ell_K & = & \text{the \, length \, of \, the \, longest \, chain \, of \, faces \, in } K, \\
\ell   & = & \min \{ \ell_K - 1, \dim \, L \}.
\end{array}
\end{equation}
We prove the correctness of our FRA and an 
upper bound on the number of reducing iterations in Theorem \ref{fra-thm}: 
\bth
\label{fra-thm}
Suppose that the FRA finds $y_0, y_1, \dots, $ and corresponding faces 
$F_0, \, F_1, \dots$ 
Then the following hold:
\begin{enumerate}
\item \label{fra-thm-1} $F_{\min} \subseteq F_i \, $ for $i=0, 1, \dots$ 
\item \label{fra-thm-2} After a sufficiently large number of reducing iterations the algorithm finds 
\mbox{$F_{\min} = F_t$} in some iteration $t. \,$ Furthermore, 
$$
F_{\min} = F_i
$$
holds for all $i \geq t.$ 
\item \label{fra-thm-3} The number of reducing iterations in the FRA is at most $\ell.$ 
\end{enumerate}
\enth
\pf{}
Let us first note that the face $F_i$ found by the algorithm is of the form
$$
F_i \, = \, K \cap y_0^\perp \cap \dots \cap y_i^\perp, \, i=0, 1, \dots
$$
Statement (\ref{fra-thm-1})  follows from applying repeatedly part 
(\ref{red-lemma-1}) of Lemma \ref{red-lemma}.

In (\ref{fra-thm-2}) the first part of the claim is straightforward; in particular, 
the number of reducing iterations cannot exceed $\ell_K -1.$ 
Suppose $i \geq t.$ Since $F_t = F_{\min}, \,$ we have 
\begin{equation} \label{ftfi}
F_{\min} \subseteq F_i = F_t \cap y_{t+1}^\perp \cap \dots \cap y_i^\perp \subseteq F_{\min}, 
\end{equation}
so equality holds throughout in (\ref{ftfi}), which proves $F_i = F_{\min}.$ 

To prove (\ref{fra-thm-3}) let us denote by $k$ the number of reducing iterations.
It remains to show that $k \leq \dim L$ holds, so assume to the contrary $k > \dim L.$
Suppose that $y_{i_1}, \dots, y_{i_k}$ are the vectors found in reducing iterations, where $i_1 < \dots < i_k.$
Since they are all in $L, \,$ they must be linearly dependent, 
so there is an index  $r \in \{1, \dots, k \}$ such that 
$$y_{i_r} \in \lin \{ y_{i_1}, \dots, y_{i_{r-1}} \} \, \subseteq \, \lin \{ y_0, y_{1}, \dots, y_{i_{r}-1} \}.$$ 
For brevity let us write $s=i_r.$ Then 
$y_0^\perp \cap \dots \cap y_{s-1}^\perp \subseteq y_s^\perp, \,$ so
$$F_{s} = F_{s-1}, \,$$ 
i.e., the $s^{th}$ step is not reducing, which is a contradiction. 
\qed 

Next we illustrate our algorithm on the examples of Section \ref{sect-prelim}.

{\bf Examples \ref{lin-ex} and \ref{ex1} continued} 
Suppose we run our algorithm on the linear system (\ref{linsystem}).
The $y_i$ vectors below, with corresponding faces shown, are a possible output: 
\renewcommand{\arraystretch}{1.1}
\begin{equation} \label{y0y1}
\begin{array}{rcl} 
y_0 & = & 0,  \,   F_0 = \radp{5},    \\
y_1 & = & (0,  \, 0,  \, 0,  \, 1,  \, 1 )^T, \, F_1 = \radp{3} \times \{ 0 \}^2, \\
y_2 & = & (0,  \, 1,  \, 1,  \, 0,  \, -1 )^T, \, F_2 = F_{\rm min} = \radp{1} \times \{ 0 \}^4.
\end{array}
\end{equation}
\renewcommand{\arraystretch}{1.0}
The algorithm may also finish in one step, by finding, say, $y_0 = 0, \,$ and 
\begin{equation}
y_1 \, = \, (0,  \, 1,  \, 1,  \, 2,  \, 1 )^T.
\end{equation}
Of course, in linear systems 
there is always a reducing certificate that finds the minimal cone in one step, 
i.e., $F_{\min} = K \cap y_1^\perp$ for some $y_1 \geq 0;$ this is straightforward from LP duality.

When we run our algorithm on the instance of 
(\ref{ex1-p}), the $y_i$ matrices below, with corresponding $F_i$ faces, 
are a possible output:
\renewcommand{\arraystretch}{1.1}
\begin{equation} \label{y0y1y2}
\begin{array}{rcl} 
y_0 & = & 0,  \,   F_0 = \psd{3},    \\
y_1 & = & \begin{pmatrix} 0 & 0 & 0 \\ 0 & 0 & 0 \\ 0 & 0 & 1 \end{pmatrix}, \, F_1 \, = \, 
\begin{pmatrix} 
  \parbox{1cm}{$\,\, \bigoplus$}  & \hspace{-.5cm} \parbox{1cm} {$\begin{array}{c} 0 \\ 0 \end{array}$} \hspace{-.5cm} \\
  \parbox{1cm}{$0 \,\,\,\, 0 \,\,\,\,$}  &   \hspace{-.45cm} 0 
\end{pmatrix}, \,\, \\
y_2 & = & \begin{pmatrix} 0 & 0 & -1 \\ 0 & 2 & 0 \\ -1 & 0 & 0 \end{pmatrix}, \, F_2 \, = \, F_{\min} \, = \, 
\begin{pmatrix} \oplus & 0 & 0 \\ 0 & 0 & 0 \\ 0 & 0 & 0 \end{pmatrix}.
\end{array}
\end{equation}
\renewcommand{\arraystretch}{1}
Indeed it is clear that the $y_i$ are orthogonal to all the constraint matrices 
in problem (\ref{ex1-p}) and that $y_i \in F_{i-1}^*$ for $i=1,2.$ 
\qed

Let us now consider the conic system 
\begin{equation} \label{ext} \tag{\mbox{$EXT$}} \left.
\begin{array}{rcl} 
y_0      & =   & 0  \\
y_i      & \in & F_{i-1}^*, \, {\rm where}    \\
F_{i-1}  & =   & K \cap y_0^\perp \cap \dots \cap y_{i-1}^\perp, \, i=1, \dots, \ell+1  \\
y_i      & \in &  L, \, i=1, \dots, \ell      
\end{array} \right\},
\end{equation}
that we call an extended system. 

We have the following representation theorem:
\bth \label{fmin-fra}
$F_{\min}^* \, = \, \{ \, y_{\ell+1} \, | \, (y_i)_{i=0}^{\ell+1} \, \text{is feasible in} \, (\ref{ext}) \, \}.$
\qed
\enth 

Before proving Theorem \ref{fmin-fra} we make some remarks. 
First, the two different ranges for the $i$ indices 
in the constraints of (\ref{ext}) are not accidental: the sequence 
$y_0, \dots, y_\ell$ is a possible output of our FRA, iff with some $y_{\ell+1}$ it is 
feasible in (\ref{ext}), and the variable $y_{\ell+1}$ represents the dual of the minimal cone.
It also becomes clearer now why we allow nonreducing iterations in our algorithm: in the 
conic system (\ref{ext}) some $y_i$ correspond to reducing iterations, but others do not.

The extended system (\ref{ext}) is not linear, due to how the $y_i$ vectors depend on the previous 
$y_j, \,$ and in general we also don't know how to describe the duals 
of faces of $K.$ Hence the representation of  Theorem \ref{fmin-fra} is not yet immediately useful.
However, in the next section we state an equivalent conic linear 
system to represent $F_{\min}^* \,$ 
when $K$ is nice, and 
arrive at the representation of (\ref{fmin*repr}), and at an extended dual of (\ref{p}).

\pf{of Theorem \ref{fmin-fra}} 
Let us write $G$ for the set on the right hand side.
Suppose that $(y_i)_{i=0}^{\ell+1}$ is feasible in (\ref{ext}) with corresponding faces 
$F_0, \dots, F_\ell.$ By part (\ref{fra-thm-1}) in Theorem \ref{fra-thm} we have 
\begin{equation} \label{bla}
F_{\min} \subseteq F_\ell, \, \text{hence} \, F_{\min}^* \supseteq F_\ell^*.
\end{equation}
Since $y_{\ell+1} \in F_\ell^* \,$ in $G, \,$ 
the containment $F_{\min}^* \supseteq G$ follows. 

By part (\ref{fra-thm-2})-(\ref{fra-thm-3}) in Theorem \ref{fra-thm} 
there exists $(y_i)_{i=0}^{\ell+1}$ that is feasible in (\ref{ext}), with corresponding faces 
$F_0, \dots, F_{\ell}$ such that equality 
holds in (\ref{bla}).
This proves the inclusion $F_{\min}^* \subseteq G.$ 
\qed

\section{When $K$ is nice: an extended dual, and an explicit extended dual for semidefinite programs}
\label{section-nice}

From now on we make the following assumption:
\begin{center}
\framebox{ $K$ is nice.}
\end{center}
Let us recall the definition of $L$ from (\ref{def-L}), 
and consider the conic system
\begin{equation} \label{extnice} \tag{\mbox{$EXT_{\rm nice}$}} \left.
\begin{array}{rcl} 
(u_0, v_0)       & =   & (0,0)  \\
(u_{i},v_{i})    & \in & K^* \times \tan(u_0+\dots+u_{i-1}, K^*), \, i=1, \dots, \ell+1 \\
u_{i}+v_{i}      & \in &  L, \, i=1, \dots, \ell
\end{array} \right\}.
\end{equation}
This is a conic linear system, since the set 
\begin{align*} 
\{ \, (u,v) \, | \, u \in K^*, \, v \in \tan(u, K^*) \, \}
\end{align*}
is a convex cone, although it may not be closed 
(e.g., if $K^* = \radp{2}, \,$ 
then $(\epsilon, 1)$ is in this set for all $\epsilon > 0, \,$ but 
$(0,1)$ is not). 
\bth 
\label{yuv} 
$
\Feas(\ref{ext}) \, = \, \{ \, (u_i + v_i)_{i=0}^{\ell+1} : (u_i, v_i)_{i=0}^{\ell+1} \, \in \, \Feas(\ref{extnice}) \, \}.
$ 
\enth 
\pf{of $\subseteq$} Suppose that $(y_i)_{i=0}^{\ell+1}$  
is feasible in (\ref{ext}), with faces 
\begin{eqnarray} \label{defFi}
F_{i-1} & = & K \cap y_0^\perp \cap \dots \cap y_{i-1}^\perp, \, i=1, \dots, \ell+1.
\end{eqnarray}
For $i=1, \dots, \ell+1$ we have $y_i \in F_{i-1}^*, $ and $K$ is nice, 
so we can write 
$y_i  =  u_i + v_i$ for some  $u_i \in K^*$ and $v_i \in F_{i-1}^\perp.$ 
Also, let us set $u_0 = v_0 = 0,$ then of course $y_0 = u_0 + v_0.$ 

We show that $(u_i, v_i)_{i=0}^{\ell+1}$ is feasible in 
(\ref{extnice}). To do this, it is enough to verify 
\begin{eqnarray}
\label{FiperpK} 
F_{i-1}^\perp & = & \tan(u_0 + \dots + u_{i-1}, K^*) \, 
\end{eqnarray}
for $i=1, \dots, \ell+1.$ 
Equation (\ref{FiperpK}) will follow if we prove 
\begin{eqnarray} \label{FiK}
F_{i-1} & = & K \cap (u_0 + \dots + u_{i-1})^\perp 
\end{eqnarray}
for $i=1, \dots, \ell+1;$ indeed, from (\ref{FiK}) we directly obtain
\begin{eqnarray*}
F_{i-1}       & = & \face(u_0 + \dots + u_{i-1}, K^*)^\triangle,
\end{eqnarray*}
hence 
\begin{eqnarray*}
F_{i-1}^\perp       & = & \face(u_0 + \dots + u_{i-1}, K^*)^{\triangle \perp} \\
                    & = & \tan(u_0 + \dots + u_{i-1}, K^*),
\end{eqnarray*}
where the second equality comes from (\ref{tanxK}). 

So it remains to prove (\ref{FiK}). It is clearly true for $i=1.$ 
Let $i$ be a nonnegative integer at most $\ell+1$ and assume that (\ref{FiK}) holds for 
$1, \dots, i-1.$ We then have 
$$
\begin{array}{rcl}
F_{i-1} & = & F_{i-2} \cap y_{i-1}^\perp \\
      & = & F_{i-2} \cap (u_{i-1} + v_{i-1})^\perp  \\
      & = & F_{i-2} \cap u_{i-1}^\perp  \\
      & = & K \cap (u_0 + \dots + u_{i-2})^\perp \cap u_{i-1}^\perp \\
      & = & K \cap (u_0 + \dots + u_{i-2} + u_{i-1})^\perp.
\end{array}
$$
Here the second equation follows from the definition of $(u_{i-1}, v_{i-1}), \,$ 
the third from $v_{i-1} \in F_{i-2}^\perp, \,$ 
the fourth from the inductive hypothesis, and the last from 
all $u_j$ being in $K^*.$ 

Thus the proof of the containment $\subseteq$ is complete.

\pf{of $\supseteq$} Let us choose $(u_i, v_i)_{i=0}^{\ell+1}$ to be feasible in (\ref{extnice}), 
define \mbox{$y_i = u_i + v_i$} for all $i, \,$ and the faces $F_0, \dots, F_\ell$ as in 
(\ref{defFi}). 
Repeating the previous argument verbatim, 
(\ref{FiperpK}) holds, so we have 
$$
y_i \in K^* + F_{i-1}^\perp = F_{i-1}^*, \, i=1, \dots, \ell+1.
$$ 
Therefore $(y_i)_{i=0}^{\ell+1}$ is 
feasible in (\ref{ext}) and this completes the proof.
\qed 

We now arrive at the representation of $F_{\min}^*$ that we previewed in (\ref{fmin*repr}), and at 
an extended dual of (\ref{p}):
\bcor \label{fmin*cor}
The dual of the minimal cone of (\ref{p}) has a representation
\begin{equation} \label{fmin*-cor-eq}
\begin{array}{rcl}
F_{\min}^* & = & \{ \, u_{\ell+1} + v_{\ell+1} \, : \,  (u_i, v_i)_{i=0}^{\ell+1} \, \text{is feasible in} \, (\ref{extnice}) \, \},
\end{array}
\end{equation}
and the extended dual 
\begin{equation}  \tag{$D_{\rm ext}$} \label{d-ext}
\begin{array}{rl} 
 \inf        &   \langle b, u_{\ell+1} + v_{\ell+1}  \rangle  \\                   
s.t.         & \, A^*(u_{\ell+1} + v_{\ell+1}) = c \\
             &  \, (u_i, v_i)_{i=0}^{\ell+1} \, \text{is feasible in} \, (\ref{extnice})
\end{array}
\end{equation}
has strong duality with (\ref{p}). 

In particular, if (\ref{p}) is a semidefinite program with $m$ variables, 
independent 
constraint matrices, and $K = \psd{n},\,$ then the problem 
\begin{equation}  \tag{$D_{\rm ext, SDP}$} \label{d-ext-sdp}
\begin{array}{rrcl} 
 \inf        &   \langle b, u_{\ell+1} + v_{\ell+1}  \rangle  \\                   
  s.t.         & \, A^*(u_{\ell+1} + v_{\ell+1}) & = & c \\
             &  \, (A, b)^*(u_i + v_i) & = & 0, \, i=1, \dots, \ell \\
             &  \, u_i & \succeq & 0, \, i=1, \dots, \ell+1 \\
             & (*) \, \begin{pmatrix} u_0 + \dots + u_{i-1} & w_i \\ w_i^T  & \beta_i I \end{pmatrix} & \succeq & 0, i=1, \dots, \ell+1 \\
             & v_i & = & w_i + w_i^T, \, i=1, \dots, \ell+1 \\
             &  \, w_i & \in & \rad{n \times n},  \, i  =  1, \dots, \ell+1 \\
             &  \, \beta_i & \in & \rad{},  \, i  =  1, \dots, \ell+1 \\
             &  \,  (u_0, v_0) & = & (0, 0),
\end{array}
\end{equation}
where 
\begin{equation} \label{def-ell-sdp}
\ell = \min \, \{ \, n, n(n+1)/2 - m -1 \, \}, 
\end{equation}
has strong duality with (\ref{p}).
\ecor
\pf{}
The representation (\ref{fmin*-cor-eq}) follows from combining 
Theorems \ref{fmin-fra} and \ref{yuv}. The second statement of the theorem follows, since 
replacing $K^*$ by $F_{\min}^*$ in (\ref{d}) yields a strong dual for (\ref{p}).

Suppose now that (\ref{p}) is a semidefinite program with $K = \psd{n}, \,$ 
with $m$ variables, and with independent constraint matrices.
The length of the longest chain of faces in $\psd{n}$ is $n+1 \,$ and 
the dimension of the subspace ${\cal N}( (A, b)^*)$ is 
$n(n+1)/2 - m -1. \,$ Hence we can choose $\ell$ as in (\ref{def-ell-sdp}) to obtain a correct
extended dual.

Let $v_i \in \sym{n}$ and $u_0, \dots, u_{i-1} \in \psd{n}, \,$ where $i \in \{ \, 1, \dots, \ell+1 \, \}.$ 
The representation of the tangent space in $\psd{n}$ 
in (\ref{tan-x}) implies that $v_i \in \tan(u_0 + \dots + u_{i-1}, K^*) \,$ holds, 
iff $v_i, u_0, \dots, u_{i-1}$ 
with some $w_i$ (possibly nonsymmetric) matrices and $\beta_i$ scalars
satisfies the $i$th constraint of (\ref{d-ext-sdp}) marked by (*). 
This proves the correctness of the extended dual (\ref{d-ext-sdp}).
\qed

For the reducing certificates found for the linear system (\ref{linsystem}) 
and displayed in (\ref{y0y1}) the reader can easily find the decomposition whose existence 
we showed in Theorem \ref{yuv}.

\noindent{\bf Example \ref{ex1} continued} Recall that when we run our FRA on the SDP instance
(\ref{ex1-p}), matrices $y_0, y_1, y_2$ shown in equation (\ref{y0y1y2}) are a possible output.

We illustrate their decomposition as proved in Theorem \ref{yuv}, 
in particular, as $y_i = u_i + v_i$ with $u_i \in K^*$ and 
$v_i \in \tan(u_0 + \dots + u_{i-1}, K^*)$ for $i=1,2: \,$ 
\renewcommand{\arraystretch}{1.1}
\begin{equation} \label{uiviex1}
\begin{array}{rcl} 
u_0 & = & 0, \, v_0 \, = \, 0,    \\
u_1 & = & \begin{pmatrix} 0 & 0 & 0 \\ 0 & 0 & 0 \\ 0 & 0 & 1 \end{pmatrix},  \, v_1 \, = \, 0, \\
u_2 & = & \begin{pmatrix} 0 & 0 & 0 \\ 0 & 2 & 0 \\ 0 & 0 & 0 \end{pmatrix}, \, v_2 \, = \, \begin{pmatrix} 0 & 0 & -1 \\ 0 & 0 & 0 \\ -1 & 0 & 0 \end{pmatrix}.
\end{array}
\end{equation}
\renewcommand{\arraystretch}{1}
We can check $v_2 \in \tan(u_1, \psd{3})$ by 
using the tangent space formula (\ref{tanir}). 

To illustrate the correctness of the extended dual (\ref{d-ext-sdp}), we first note 
that 
$n=m=3, \,$ so by formula (\ref{def-ell-sdp}) we can choose $\ell=2$ to obtain a correct 
extended dual. 
Recall that $y \in F_{\min}^*$ is an optimal dual solution 
if and only if it is of the form 
\begin{equation} \label{yfmin}
\begin{array}{rcl}
y & = &  \begin{pmatrix} 0  & 1/2       & - y_{22}/2 \\
           1/2    & y_{22}   & y_{23} \\
     - y_{22}/2   &  y_{23}  & y_{33} \end{pmatrix}.
\end{array}
\end{equation}
Consider the $(u_i, v_i)_{i=0}^{2}$ sequence shown in (\ref{uiviex1});
we prove  that any $y$ optimal matrix satisfies 
\begin{equation} \label{jerry}
y \, \in \, \psd{3} + \tan(u_0 + u_1 + u_2, \psd{3}).
\end{equation}
Indeed, 
$\tan(u_0 + u_1 + u_2, \psd{3})$ is the set of $3$ by $3$ matrices with the component in the 
$(1,1)$ position equal to zero, and the other components arbitrary, and this proves (\ref{jerry}).

In fact, considering the expression for $F_{\min}^*$ in (\ref{fminex1}), 
it follows that any $y \in F_{\min}^*$ can be decomposed as in equation (\ref{jerry}).
\qed

\section{Variants of extended duals } 
\label{section-variants}
So far we proved the correctness of an extended dual of (\ref{p}), which is itself an explicit 
semidefinite program when (\ref{p}) is. Ramana's original dual is somewhat different from
(\ref{d-ext-sdp}) though. Here we describe several variants of extended duals for 
(\ref{p}) and show how to derive Ramana's dual.

First let us define a simplified extended system
\begin{equation} \label{extnicep} \tag{\mbox{$EXT_{\rm nice, simple}$}} \left.
\begin{array}{rcl} 
(u_0, v_0)       & =   & (0,0)  \\
(u_{i},v_{i})    & \in & K^* \times \tan(u_{i-1}, K^*), \, i=1, \dots, \ell+1 \\
u_{i}+v_{i}      & \in &  L, \, i=1, \dots, \ell.      
\end{array} \right\}
\end{equation}
We prove that this system works just as well as (\ref{extnice}) when constructing extended duals.
\bcor
\label{fmin*cor-simple}
The dual of the minimal cone of (\ref{p}) has a representation
\begin{equation} \label{fmin*-cor-eq-sh}
\begin{array}{rcl}
F_{\min}^* & = & \{ \, u_{\ell+1} + v_{\ell+1} \, : \,  (u_i, v_i)_{i=0}^{\ell+1} \, \text{is feasible in} \, (\ref{extnicep}) \, \},
\end{array}
\end{equation}
and the extended dual 
\begin{equation}  \tag{$D_{\rm ext, simple}$} \label{d-ext-p}
\begin{array}{rl} 
 \inf        &   \langle b, u_{\ell+1} + v_{\ell+1}  \rangle  \\                   
s.t.         & \, A^*(u_{\ell+1} + v_{\ell+1}) = c \\
             &  \, (u_i, v_i)_{i=0}^{\ell+1} \, \text{is feasible in} \, (\ref{extnicep}), 
\end{array}
\end{equation}
where $\ell$ is defined in (\ref{def-ell}), has strong duality with (\ref{p}). 

In particular, if (\ref{p}) is an SDP as described in Corollary \ref{fmin*cor}, 
then the problem obtained from 
(\ref{d-ext-sdp}) by replacing the constraint (*) by 
$$
(**) \, \begin{pmatrix} u_{i-1} & w_i \\ w_i^T  & \beta_i I \end{pmatrix} \, \succeq \, 0, \; i=1, \dots, \ell+1,
$$
has strong duality with (\ref{p}). 
\ecor
\pf{} It is enough to prove the representation in equation 
(\ref{fmin*-cor-eq-sh}); given this, the rest of the 
proof is  analogous to the proof of the second and third statements in Corollary
\ref{fmin*cor}. 

We will use the representation of $F_{\min}^*$ in (\ref{fmin*-cor-eq}).
Let us denote by $G$ the set on the right hand side of equation 
(\ref{fmin*-cor-eq-sh}); we will prove $G = F_{\min}^*.$

To show  $G \subseteq F_{\rm min}^*$ suppose 
$u_{\ell+1} + v_{\ell+1} \in G, $ where 
$(u_i, v_i)_{i=0}^{\ell+1}$ is feasible in (\ref{extnicep}). 
Then it is also feasible in (\ref{extnice}), since applying the tangent space formula 
(\ref{tanxK}) with $K^*$ in place of $K$ implies that 
$$\tan(u_{i-1}, K^*) \, \subseteq \, \tan(u_0 + \dots + u_{i-1}, K^*) \,$$ 
holds for $i=1, \dots, \ell+1.$ 

To prove $G \supseteq F_{\rm min}^*$ suppose that 
$u_{\ell+1} + v_{\ell+1} \in F_{\rm min}^*, $ where 
$(u_i, v_i)_{i=0}^{\ell+1}$ is feasible in (\ref{extnice}). 
Again, by (\ref{tanxK}) the sets $\tan(u_0, K^*), \dots,$ $\tan(u_0 + \dots + u_{i-1}, K^*)$ 
are all contained 
in $\tan(u_0 + \dots + u_{i-1}, K^*) \,$ for $i=1, \dots, \ell.$
Hence 
\begin{equation} \label{v1dots} 
v_1 + \dots + v_i \, \in \, \tan(u_0 + \dots + u_{i-1}, K^*),  \, i=1, \dots, \ell \\
\end{equation}
holds, 
and we also have
\begin{equation} \label{prakash}
v_{\ell+1} \, \in \, \tan(u_0 + \dots + u_{\ell}, K^*).
\end{equation}
Let us define 
$$
(u_i', v_i')  \, = \, (u_0 + \dots + u_i, v_0 + \dots + v_i), \, i=1, \dots, \ell.
$$
By (\ref{v1dots}) and (\ref{prakash}) it follows that 
$(u_{\ell+1}, v_{\ell+1})$ with $(u_i', v_i')_{i=0}^{\ell}$ is feasible for (\ref{extnicep}),
so the inclusion follows.
\qed

Let us now consider another extended system 
\begin{equation} \label{extnicepttan} \tag{\mbox{$EXT'_{\rm nice, simple}$}} \left.
\begin{array}{rcl} 
(u_0, v_0)       & =   & (0,0)  \\
(u_{i},v_{i})    & \in & K^* \times \ttan(u_{i-1}, K^*), \, i=1, \dots, \ell+1 \\
u_{i}+v_{i}      & \in &  L, \, i=1, \dots, \ell
\end{array} \right\},
\end{equation}
where the set $\ttan(u, K^*), $ 
satisfies the following two requirements for all $u \in K^*:$ 
\begin{enumerate}
\item $\ttan(u, K^*) \subseteq \tan(u, K^*).$
\item For all $v \in \tan(u, K^*)$ there exists $\lambda_v >  0$ such that $v \in \ttan(\lambda_v u, K^*).$ 
\end{enumerate}
\bcor
The dual of the minimal cone of (\ref{p}) has the representation
\begin{equation} \label{fmin*-cor-eq-ttan}
\begin{array}{rcl}
F_{\min}^* & = & \{ \, u_{\ell+1} + v_{\ell+1} \, : \,  (u_i, v_i)_{i=0}^{\ell+1} \, \text{is feasible in} \,\, (\ref{extnicepttan}) \, \},
\end{array}
\end{equation}
and the extended dual 
\begin{equation}  \tag{$D'_{\rm ext, simple}$} \label{d-ext-p-ttan}
\begin{array}{rl} 
 \inf        &   \langle b, u_{\ell+1} + v_{\ell+1}  \rangle  \\                   
s.t.         & \, A^*(u_{\ell+1} + v_{\ell+1}) = c \\
             &  \, (u_i, v_i)_{i=0}^{\ell+1} \, \text{is feasible in} \, (\ref{extnicepttan}),
\end{array}
\end{equation}
where $\ell$ is defined in (\ref{def-ell}), has strong duality with (\ref{p}). 

In particular, if (\ref{p}) is an SDP as described in Corollary \ref{fmin*cor}, 
then the problem obtained from 
(\ref{d-ext-sdp}) by replacing the constraint (*) by 
$$
(***) \, \begin{pmatrix} u_{i-1} & w_i \\ w_i^T  & I \end{pmatrix} \, \succeq \; 0, \; i=1, \dots, \ell+1,
$$
and dropping the $\beta_i$ variables, has strong duality with (\ref{p}). 
\ecor

\pf{}  We use the representation of $F_{\min}^*$ in (\ref{fmin*-cor-eq-sh}).
Let us denote by $G$ the set on the right hand side of equation 
(\ref{fmin*-cor-eq-ttan}). We will prove $G = F_{\min}^*.$

It is clear that $G \subseteq F_{\min}^*,$ since if 
$(u_i, v_i)_{i=0}^{\ell+1}$ is feasible in (\ref{extnicepttan}), then 
by the first property of the operator  tan' it is also feasible in 
(\ref{extnicep}).

To show the opposite inclusion, suppose 
$u_{\ell+1} + v_{\ell+1} \in F_{\min}^*, $ where 
$(u_i, v_i)_{i=0}^{\ell+1}$ is feasible in (\ref{extnicep}). 
Let us choose $\lambda_{\ell}, \lambda_{\ell-1}, \dots, \lambda_1$ positive reals 
such that 
\begin{equation}
\begin{array}{rcl}
v_{\ell+1}              & \in & \ttan(\lambda_{\ell} u_{\ell}, K^*), \\
\lambda_{\ell} v_{\ell} & \in & \ttan(\lambda_{\ell-1} u_{\ell-1}, K^*), \\
                        & \vdots &                                      \\
\lambda_{2} v_{2}       & \in & \ttan(\lambda_{1} u_{1}, K^*),
\end{array}
\end{equation}
and for completeness, set $\lambda_0 = 0.$ 
Then $(u_{\ell+1}, v_{\ell+1})$ with 
$(\lambda_i u_i, \lambda_i v_i)_{i=0}^{\ell}$ is feasible in 
(\ref{extnicepttan}), and this proves $F_{\min}^* \subseteq G.$ 
\qed

We finally remark that in the extended duals for semidefinite programming
it is possible to eliminate the $v_i$ variables and use the $w_i$ matrices 
directly in the constraints; thus one can exactly obtain Ramana's dual.
We leave the details to the reader.

\section{Conclusion}
\label{repr-sect}

We gave a simple and self-contained exposition of a facial reduction algorithm and of extended duals: 
both approaches yield strong duality for a conic linear program, without assuming any constraint 
qualification. We generalized Ramana's dual: 
we proved that when $K$ is a nice cone, the set $F_{\min}^*$ 
has an extended formulation, i.e., it is 
the projection of the feasible set of a conic linear system in a higher dimensional space.
The only nontrivial constraints in this system
are of the form $u \geq_{K^*} 0$, and $v \in \tan(u, K^*)$.

This formulation leads to  an extended, strong  dual of (\ref{P}), when $K$ is nice.
When $K = K^*$ is the semidefinite cone, by writing the tangent space constraint 
as a semidefinite constraint, we obtain an extended strong dual, which is an SDP itself, and thus 
recover variants of Ramana's dual.

One may wonder, whether $F_{\min}$ itself has an extended formulation.
Suppose that $K$ is an arbitary closed convex cone.
When a fixed $\bar{s} \in \ri F_{\min}$ is given, then obviously
\begin{equation*}
F_{\min} \, = \, \{ \, s \, | \, 0  \leq_K   s  \leq_K  \alpha \bar{s}  \;  \text{for  some} \; \alpha \geq 0 \, \}.
\end{equation*}
We can also represent the minimal cone without such an $\bar{s}$, since 
\begin{equation} \label{fmin2}
F_{\min} \, = \, \{ \, s \, | \, 0  \leq_K   s  \leq_K  \alpha b - Ax   \;  \text{for  some} \; x, \; \text{and} \; \alpha \geq 0 \, \}.
\end{equation}
This representation was obtained by Freund \cite{Freund:94}, 
based on the article by himself, Roundy and Todd \cite{FreundRoundyTodd:85}.

It is also natural to ask, whether there are other nice cones, 
for which the set
\begin{align*} 
\{ \, (u,v) \, | \, u \in K^*, \, v \in \tan(u, K^*) \, \}
\end{align*}
has a formulation in terms of $K^*;$ e.g., 
is this true for the second order cone? 
Conic linear programs over such cones would 
also have Ramana-type (ie., expressed only in terms of $K^*$) extended duals.

\noindent{\bf Acknowledgement} I would like to thank an anonymous referee and Minghui Liu 
for their helpful comments.

\bibliography{D:/bibfiles/mysdp}

\begin{thebibliography}{10}

\bibitem{AliGold:03}
Farid Alizadeh and Donald Goldfarb.
\newblock Second-order cone programming.
\newblock {\em Math. Program. Ser. B}, 95:3--51, 2003.

\bibitem{AndRoosTerlaky:02}
Erling~D. Andersen, Cees Roos, and Tam\'as Terlaky.
\newblock Notes on duality in second order and {$p$}-order cone optimization.
\newblock {\em Optimization}, 51(4):627--643, 2002.

\bibitem{BentalNem:01}
Aharon Ben-Tal and Arkadii Nemirovskii.
\newblock {\em Lectures on modern convex optimization}.
\newblock MPS/SIAM Series on Optimization. SIAM, Philadelphia, PA, 2001.

\bibitem{BonnShap:00}
Fr\'ed\'eric~J. Bonnans and Alexander Shapiro.
\newblock {\em Perturbation analysis of optimization problems}.
\newblock Springer Series in Operations Research. Springer-Verlag, 2000.

\bibitem{BorWolk:81B}
Jonathan~M. Borwein and Henry Wolkowicz.
\newblock Facial reduction for a cone-convex programming problem.
\newblock {\em J. Aust. Math. Soc.}, 30:369--380, 1981.

\bibitem{BorWolk:81}
Jonathan~M. Borwein and Henry Wolkowicz.
\newblock Regularizing the abstract convex program.
\newblock {\em J. Math. Anal. App.}, 83:495--530, 1981.

\bibitem{BoydVand:04}
Stephen Boyd and Lieven Vandenberghe.
\newblock {\em Convex Optimization}.
\newblock Cambridge University Press, 2004.

\bibitem{CheWolkSchurr:12}
Vris Cheung, Henry Wolkowicz, and Simon Schurr.
\newblock Preprocessing and regularization for degenerate semidefinite
  programs.
\newblock Technical report, University of Waterloo, 2012.

\bibitem{Chua:03}
Check-Beng Chua.
\newblock Relating homogeneous cones and positive definite cones via
  {T}-algebras.
\newblock {\em SIAM J. Optim.}, 14:500--506, 2003.

\bibitem{ChuaTuncel:08}
Check-Beng Chua and Levent Tun\c{c}el.
\newblock Invariance and efficiency of convex representations.
\newblock {\em Math. Program. B}, 111:113--140, 2008.

\bibitem{Faybu:02}
Leonid Faybusovich.
\newblock On {N}esterov's approach to semi-definite programming.
\newblock {\em Acta Appl. Math.}, 74:195--215, 2002.

\bibitem{Freund:94}
Robert~M. Freund.
\newblock Talk at the {U}niversity of {W}aterloo.
\newblock 1994.

\bibitem{FreundRoundyTodd:85}
Robert~M. Freund, Robin Roundy, and Michael~J. Todd.
\newblock Identifying the set of always-active constraints in a system of
  linear inequalities by a single linear program.
\newblock Technical Report Sloan W.P. No. 1674-85, Sloan School of Business,
  MIT, 1985.

\bibitem{Guler:10}
Osman G$\ddot{\mathrm{u}}$ler.
\newblock {\em Foundations of Optimization}.
\newblock Graduate Texts in Mathematics. Springer, 2010.

\bibitem{Glineur:01}
Francois Glineur.
\newblock Proving strong duality for geometric optimization using a conic
  formulation.
\newblock {\em Ann. Oper. Res.}, 105(2):155--184, 2001.

\bibitem{KlepSchw:12}
Igor Klep and Markus Schweighofer.
\newblock An exact duality theory for semidefinite programming based on sums of
  squares.
\newblock Technical report, Universit{\"{a}}t Konstanz, 2012,
  http://arxiv.org/abs/1207.1691.

\bibitem{Lewis:94}
Adrian~S. Lewis.
\newblock Facial reduction in partially finite convex programming.
\newblock {\em Math. Program. B}, 65:123--138, 1994.

\bibitem{LuoSturm:00}
Zhi-Quan Luo and Jos Sturm.
\newblock Error analysis.
\newblock In Romesh Saigal, Lieven Vandenberghe, and Henry Wolkowicz, editors,
  {\em Handbook of semidefinite programming}. Kluwer Academic Publishers, 2000.

\bibitem{LuoSturmZhang:97}
Zhi-Quan Luo, Jos Sturm, and Shuzhong Zhang.
\newblock Duality results for conic convex programming.
\newblock Technical Report Report 9719/A, Erasmus University Rotterdam,
  Econometric Institute, The Netherlands, 1997.

\bibitem{GouveiaParriloThomas:12}
\mbox{Jo\~{a}o} Gouveia, Pablo Parrilo, and Rekha Thomas.
\newblock Lifts of convex sets and cone factorizations.
\newblock {\em Math. Oper. Res.}, to appear.

\bibitem{Pataki:00B}
G\'abor Pataki.
\newblock A simple derivation of a facial reduction algorithm and extended dual
  systems.
\newblock Technical report, Columbia University, 2000.

\bibitem{Pataki:07}
G{\'a}bor Pataki.
\newblock On the closedness of the linear image of a closed convex cone.
\newblock {\em Math. Oper. Res.}, 32(2):395--412, 2007.

\bibitem{Pataki:10}
G\'abor Pataki.
\newblock Bad semidefinite programs: they all look the same.
\newblock Technical Report available from http://arxiv.org/abs/1112.1436,
  University of North Carolina at Chapel Hill, 2010.

\bibitem{Pataki:12}
G\'abor Pataki.
\newblock On the connection of facially exposed and nice cones.
\newblock {\em J. Math. Anal. App.}, 400:211--221, 2013.

\bibitem{PolikTerlaky:09}
Imre P\'olik and Tam\'as Terlaky.
\newblock Exact duality for optimization over symmetric cones.
\newblock Technical report, Lehigh University, Betlehem, PA, USA, 2009.

\bibitem{Ramana:97}
Motakuri~V. Ramana.
\newblock An exact duality theory for semidefinite programming and its
  complexity implications.
\newblock {\em Math. Program. Ser. B}, 77:129--162, 1997.

\bibitem{RaTuWo:97}
Motakuri~V. Ramana, Levent Tun\c{c}el, and Henry Wolkowicz.
\newblock Strong duality for semidefinite programming.
\newblock {\em SIAM J. Opt.}, 7(3):641--662, 1997.

\bibitem{Rockafellar:70}
Tyrrel~R. Rockafellar.
\newblock {\em Convex Analysis}.
\newblock Princeton University Press, Princeton, NJ, USA, 1970.

\bibitem{Vera:13}
Vera Roshchina.
\newblock Facially exposed cones are not nice in general.
\newblock Technical report, University of Ballarat,
  http://arxiv.org/abs/1301.1000.

\bibitem{SaigVandWolk:00}
Romesh Saigal, Lieven Vandenberghe, and Henry Wolkowicz, editors.
\newblock {\em Handbook of semidefinite programming}.
\newblock Kluwer Academic Publishers, 2000.

\bibitem{Schurretal:07}
Simon~P. Schurr, Andr\'{e}~L. Tits, and Dianne~P. O'Leary.
\newblock Universal duality in conic convex optimization.
\newblock {\em Math. Program. Ser. A}, 109:69--88, 2007.

\bibitem{Sturm:00}
Jos Sturm.
\newblock Error bounds for linear matrix inequalities.
\newblock {\em SIAM J. Optim.}, 10:1228--1248, 2000.

\bibitem{Todd:00}
Michael~J. Todd.
\newblock Semidefinite optimization.
\newblock {\em Acta Numer.}, 10:515--560, 2001.

\bibitem{Tuncel:11}
Levent Tun\c{c}el.
\newblock {\em Polyhedral and Semidefinite Programming Methods in Combinatorial
  Optimization}.
\newblock Fields Institute Monographs, 2011.

\bibitem{TunWolk:12}
Levent Tun\c{c}el and Henry Wolkowicz.
\newblock Strong duality and minimal representations for cone optimization.
\newblock {\em Comput. Optim. Appl.}, 53:619--648, 2012.

\bibitem{Waki:12}
Hayato Waki.
\newblock How to generate weakly infeasible semidefinite programs via
  lasserre's relaxations for polynomial optimization.
\newblock {\em Optim. Lett.}, 6(8):1883--1896, 2012.

\bibitem{WakiMura:12}
Hayato Waki and Masakazu Muramatsu.
\newblock Facial reduction algorithms for conic optimization problems.
\newblock {\em J. Optim. Theory Appl.}, to
  appear,DOI:10.1007/s10957-012-0219-y, 2013.

\bibitem{WangXiuLuo:11}
Yingnan Wang, Naihua Xiu, and Ziyan Luo.
\newblock A regularized strong duality for nonsymmetric semidefinite least
  squares problem.
\newblock {\em Optim. Lett.}, (5):665--682, 2011.

\end{thebibliography}

\end{document}